\def\TheMagstep{\magstep1}	
\def\PaperSize{letter}		

\magnification=\magstep1

\let\:=\colon  
   \let\?=\overline

\let\Sum=\sum \def\sum{\Sum\nolimits}

\def\IC{{\bf C}} 
\def\IP{{\bf P}}

\def\and{\hbox{ and }}
\def\Wf{\hbox{\rm W$_f$}}
	
\def\Af{\hbox{\rm A$_f$}}

\def\DONE{*!*}
\def\NextDef #1 {\def\NextOne{#1}%
 \ifx\NextOne\DONE\let\next\relax
 \else\expandafter\xdef\csname#1\endcsname{\TheOp}
  \let\next\NextDef
 \fi \next}
\def\TheOp{\mathop{\rm\NextOne}}
 \NextDef 
  Projan Supp Proj Sym Spec Hom cod Ker dist
 *!*
\def\TheOp{{\cal\NextOne}}
\NextDef 
  E F G H I J M N O R S
 *!*
\def\TheOp{\hbox{\rm\NextOne}}
\NextDef 
 A ICIS 
 *!*

\def\item#1 {\par\indent\indent\indent\indent \hangindent4\parindent
 \llap{\rm (#1)\enspace}\ignorespaces}
 \def\inpart#1 {{\rm (#1)\enspace}\ignorespaces}
 \def\part {\par\inpart}

\catcode`\@=11		

\def\vfootnote#1{\insert\footins\bgroup
 \eightpoint 
 \interlinepenalty\interfootnotelinepenalty
  \splittopskip\ht\strutbox 
  \splitmaxdepth\dp\strutbox \floatingpenalty\@MM
  \leftskip\z@skip \rightskip\z@skip \spaceskip\z@skip \xspaceskip\z@skip
  \textindent{#1}\footstrut\futurelet\next\fo@t}

\def\p.{p.\penalty\@M \thinspace}
\def\pp.{pp.\penalty\@M \thinspace}
\def\(#1){{\rm(#1)}}\let\leftp=(
\def\activeleftp{\catcode`\(=\active}
{\activeleftp\gdef({\ifmmode\let\next=\leftp \else\let\next=\(\fi\next}}

\def\sct#1\par
  {\removelastskip\vskip0pt plus2\normalbaselineskip \penalty-250 
  \vskip0pt plus-2\normalbaselineskip \bigskip
  \centerline{\smc #1}\medskip}

\newcount\sctno \sctno=0
\def\sctn{\advance\sctno by 1 
 \sct\number\sctno.\quad\ignorespaces}

\def\dno#1${\eqno\hbox{\rm(\number\sctno.#1)}$}
\def\Cs#1){\unskip~{\rm(\number\sctno.#1)}}

\def\proclaim#1 #2 {\medbreak
  {\bf#1 (\number\sctno.#2)}\enspace\bgroup\activeleftp
\it}
\def\endproclaim{\par\egroup\medskip}
\def\pf{\endproclaim{\bf Proof.}\enspace}
\def\lem{\proclaim Lemma } \def\prp{\proclaim Proposition }
\def\cor{\proclaim Corollary }	\def\thm{\proclaim Theorem }
\def\rmk#1 {\medbreak {\bf Remark (\number\sctno.#1)}\enspace}
\def\eg#1 {\medbreak {\bf Example (\number\sctno.#1)}\enspace}

\parskip=0pt plus 1.75pt \parindent10pt
\hsize29pc
\vsize44pc
\abovedisplayskip6pt plus6pt minus2pt
\belowdisplayskip6pt plus6pt minus3pt

\def\TRUE{TRUE}	
\ifx\DoublepageOutput\TRUE \def\TheMagstep{\magstep0} \fi
\mag=\TheMagstep

\newskip\vadjustskip \vadjustskip=0.5\normalbaselineskip
\def\centertext
 {\hoffset=\pgwidth \advance\hoffset-\hsize
  \advance\hoffset-2truein \divide\hoffset by 2\relax
  \voffset=\pgheight \advance\voffset-\vsize
  \advance\voffset-2truein \divide\voffset by 2\relax
  \advance\voffset\vadjustskip
 }
\newdimen\pgwidth\newdimen\pgheight
\def\letter{letter}\def\AFour{AFour}
\ifx\PaperSize\letter
 \pgwidth=8.5truein \pgheight=11truein 
 \message{- Got a paper size of letter.  }\centertext 
\fi
\ifx\PaperSize\AFour
 \pgwidth=210truemm \pgheight=297truemm 
 \message{- Got a paper size of AFour.  }\centertext
\fi

 \newdimen\fullhsize \newbox\leftcolumn
 \def\fulline{\hbox to \fullhsize}
\def\doublepageoutput
{\let\lr=L
 \output={\if L\lr
          \global\setbox\leftcolumn=\columnbox \global\let\lr=R%
        \else \doubleformat \global\let\lr=L\fi
        \ifnum\outputpenalty>-20000 \else\dosupereject\fi}%
 \def\doubleformat{\shipout\vbox{%
        \fulline{\hfil\hfil\box\leftcolumn\hfil\columnbox\hfil\hfil}%
				}%
		  }%
 \def\columnbox{\vbox
   {\makeheadline\pagebody\makefootline\advancepageno}%
   }%
 \fullhsize=\pgheight \hoffset=-1truein
 \voffset=\pgwidth \advance\voffset-\vsize
  \advance\voffset-2truein \divide\voffset by 2
  \advance\voffset\vadjustskip
 \let\firstheadline=\hfil
 
}
\ifx\DoublepageOutput\TRUE \doublepageoutput \fi

 \font\twelvebf=cmbx12		
 \font\smc=cmcsc10		

\def\eightpoint{\eightpointfonts
 \setbox\strutbox\hbox{\vrule height7\p@ depth2\p@ width\z@}%
 \eightpointparameters\eightpointfamilies
 \normalbaselines\rm
 }
\def\eightpointparameters{%
 \normalbaselineskip9\p@
 \abovedisplayskip9\p@ plus2.4\p@ minus6.2\p@
 \belowdisplayskip9\p@ plus2.4\p@ minus6.2\p@
 \abovedisplayshortskip\z@ plus2.4\p@
 \belowdisplayshortskip5.6\p@ plus2.4\p@ minus3.2\p@
 }
\newfam\smcfam
\def\eightpointfonts{%
 \font\eightrm=cmr8 \font\sixrm=cmr6
 \font\eightbf=cmbx8 \font\sixbf=cmbx6
 \font\eightit=cmti8 
 \font\eightsmc=cmcsc8
 \font\eighti=cmmi8 \font\sixi=cmmi6
 \font\eightsy=cmsy8 \font\sixsy=cmsy6
 \font\eightsl=cmsl8 \font\eighttt=cmtt8}
\def\eightpointfamilies{%
 \textfont\z@\eightrm \scriptfont\z@\sixrm  \scriptscriptfont\z@\fiverm
 \textfont\@ne\eighti \scriptfont\@ne\sixi  \scriptscriptfont\@ne\fivei
 \textfont\tw@\eightsy \scriptfont\tw@\sixsy \scriptscriptfont\tw@\fivesy
 \textfont\thr@@\tenex \scriptfont\thr@@\tenex\scriptscriptfont\thr@@\tenex
 \textfont\itfam\eightit	\def\it{\fam\itfam\eightit}%
 \textfont\slfam\eightsl	\def\sl{\fam\slfam\eightsl}%
 \textfont\ttfam\eighttt	\def\tt{\fam\ttfam\eighttt}%
 \textfont\smcfam\eightsmc	\def\smc{\fam\smcfam\eightsmc}%
 \textfont\bffam\eightbf \scriptfont\bffam\sixbf
   \scriptscriptfont\bffam\fivebf	\def\bf{\fam\bffam\eightbf}%
 \def\rm{\fam0\eightrm}%
 }

\def\today{\ifcase\month\or	
 January\or February\or March\or April\or May\or June\or
 July\or August\or September\or October\or November\or December\fi
 \space\number\day, \number\year}
\headline={%
  \ifnum\pageno=1\firstheadline
  \else
    \ifodd\pageno\oddheadline
    \else\evenheadline\fi
  \fi
}
\let\firstheadline\hfill
\def\oddheadline{\eightpoint \rlap
 \hfil\headtitle\hfil\llap{\folio}}
\def\evenheadline{\eightpoint\rlap{\folio}
 \hfil\author\hfil\llap{\today}}
\def\headtitle{\title}

 \newcount\refno \refno=0	 \def\NoKey{*!*}
 \def\MakeKey{\advance\refno by 1 \expandafter\xdef
  \csname\TheKey\endcsname{{\number\refno}}\NextKey}
 \def\NextKey#1 {\def\TheKey{#1}\ifx\TheKey\NoKey\let\next\relax
  \else\let\next\MakeKey \fi \next}
 \def\RefKeys #1\endRefKeys{\expandafter\NextKey #1 *!* }
\def\SetRef#1 #2,#3\par{%
 \hang\llap{[\csname#1\endcsname]\enspace}%
  \ignorespaces{\smc #2,}
  \ignorespaces#3\unskip.\endgraf
 }
 \newbox\keybox \setbox\keybox=\hbox{[8]\enspace}
 \newdimen\keyindent \keyindent=\wd\keybox
\def\references{
  \bgroup   \frenchspacing   \eightpoint
   \parindent=\keyindent  \parskip=\smallskipamount
   \everypar={\SetRef}}
\def\endreferences{\egroup}

 \def\serial#1#2{\expandafter\def\csname#1\endcsname ##1 ##2 ##3
  {\unskip\ #2 {\bf##1} (##2), ##3}}
 \serial{ajm}{Amer. J. Math.}
  \serial {aif} {Ann. Inst. Fourier}
 \serial{asens}{Ann. Scient. \'Ec. Norm. Sup.}
 \serial{comp}{Compositio Math.}
 \serial{conm}{Contemp. Math.}
 \serial{crasp}{C. R. Acad. Sci. Paris}
 \serial{dlnpam}{Dekker Lecture Notes in Pure and Applied Math.}
 \serial{faa}{Funct. Anal. Appl.}
 \serial{invent}{Invent. Math.}
 \serial{ma}{Math. Ann.}
 \serial{mpcps}{Math. Proc. Camb. Phil. Soc.}
 \serial{ja}{J. Algebra}
 \serial{splm}{Springer Lecture Notes in Math.}
 \serial{tams}{Trans. Amer. Math. Soc.}

\def\UThin{\penalty\@M \thinspace\ignorespaces}
\def\relaxnext@{\let\next\relax}
\def\cite#1{\relaxnext@
 \def\nextiii@##1,##2\end@{\unskip\space{\rm[\SetKey{##1},\let~=\UThin##2]}}%
 \in@,{#1}\ifin@\def\next{\nextiii@#1\end@}\else
 \def\next{{\rm[\SetKey{#1}]}}\fi\next}
\newif\ifin@
\def\in@#1#2{\def\in@@##1#1##2##3\in@@
 {\ifx\in@##2\in@false\else\in@true\fi}%
 \in@@#2#1\in@\in@@}
\def\SetKey#1{{\bf\csname#1\endcsname}}

\catcode`\@=12 

\def\title{The Multiplicity-Polar Theorem }
\def\author{Terence Gaffney}
\RefKeys BMM  E-GZ1 E-GZ2 E-GZ3 EGZS E F G2 G3  G5 G6 G-G  GK GK2 Gr H-M-S KT1 KT2 LT Mo1 MoP M    P1 P2  R T-1 T-2
 \endRefKeys

\def\topstuff{\leavevmode
 \bigskip\bigskip
 \centerline{\twelvebf \title}
 \bigskip
 \centerline{\author}
 \medskip\centerline{\today}
\bigskip\bigskip}
\topstuff

The heart of intersection theory resides in the following situation. Given a purely d-dimensional analytic set $X$, and subsets $V$ and $W$ intersecting at a point $x\in X$, $V$ and $W$ given by $i$ and $j$ equations respectively where $i+j=d$, then the correct number of times to count $x$ is the multiplicity of the ideal $I$ in the local ring of $X$ at $x$, where $I$ is the ideal generated by the equations of $V$ and $W$. The multiplicity is the correct number, because when we deform $V$ so that it intersects $W$ transversely at smooth points of $W$, the number of points is exactly the multiplicity of $I$.

In many situations it is desirable to solve the following intersection problem. Given a vector bundle of rank $e$, on a purely $d$-dimensional analytic set $X$, and two collections of sections of the bundle, such the total number of sections is $d+e-1$ and an isolated point $x\in X$, where the sections fail to have rank $e$, then how many times should we count $x$, if we are counting the number of points where our two collections fail to have maximal rank?  This problem occurs in calculating the contribution of a point to the top Chern class of a bundle on $X$ or the index of a differential 1-form with an isolated singularity on $X$, $X$ with an isolated singularity.

We may ask for more. First, just as in the case of intersecting sets we replace functions by ideals, we replace collections of vector fields by the module $M$ they generate.  Suppose the set where the rank of $M$ is less than maximal is non-isolated, but in a natural way from the problem being considered, the module $M$ is a submodule of a module $N$ where $M$ is of finite colength inside $N$.  An example of this situation is given by a 1-form on a space with isolated singularities, which is not a complete intersection, in which the 1-form still has an isolated singularity in some sense.  We would like to calculate the index of the 1-form, in a way which fits well with deformations of the original 1-form.  (In this case, the modules $M$ and $N$ are described in section 2 in the material around Theorem 2.8.)

In this paper we show that the correct solution to this problem is the multiplicity of a pair of  modules; we show that this idea is correct by proving a theorem, the multiplicity-polar theorem (Corollary 1.4) which shows that the multiplicity of a pair deforms correctly, just as the multiplicity of an ideal does.  An advantage of our theory is that it does not require the local rings  of the analytic spaces to be Cohen-Macaulay. It suffices that they be equidimensional.

If a submodule of a free module has non-finite colength, then it is tempting to generalize the Buchsbaum-Rim multiplicity of a module by using a definition based on the intersection theory approach to the multiplicity that appears in \cite{KT1}. This definition is reviewed in section 1, where we discuss the definition of the multiplicity of the pair. This generalization  yields a sequence of numbers corresponding to strata of different dimension in the co-suport of the module. This approach is described in \cite{G6}, in sections 3 and 4 and reviewed in section 1. Unfortunately, this approach does not give fiber-wise control of the integral closure of a family of modules, as is shown in section 4 of \cite{G6}.  

Expanding on this point a little, one of the goals of extending the theory of multiplicities to modules of non-finite colength is to prove a principle of specialization of integral dependence, an idea formulated in the ideal case by Bernard Teissier in \cite{T-1}. Teissier showed that if the multiplicity of a family of ideals was independent of parameter, and the restriction of a function germ to the members of the family was integrally dependent on the ideals in the family for a generic set of parameter values then the function was integrally dependent on the ideal in the local ring of the total space of the family. This result was the main tool used by Teissier in the theorems on equisingularity in \cite{T-1}. The principle of specialization of integral dependence was extended to families of modules of finite colength in \cite{GK}, and to families of ideals of non-finite colength in \cite {G-G}. The method of proof is similar. In both cases, the original space is transformed in a way that incorporates the geometry of the ideal or module. In the ideal case we blowup by the ideal, in the module case we use the $M$-transform, which is defined in section 1. If the dimension of the fiber of this transform over the origin can be controlled then the desired integral dependence result follows. The papers  \cite{GK},  \cite {G-G} both show that a jump in fiber dimension forces a jump in the invariants.  
 
The generalization of the Buchsbaum-Rim multiplicity described in \cite{G6} (which  is also generalization of the ideas of  \cite{G-G})  fails to detect the change in fiber dimension as  \cite{G6} shows.

The problem that comes up here is one of scales; there can be changes, for example, in 1-dimensional strata that affect the invariant in dimension $0$, and in fact, off set changes in the zero dimensional stratum that you want to observe, so that the zero dimensional invariant doesn't change even though the fiber dimension does.

In some problems, such as the mixed multiplicity problem for ideals of non-finite colength this problem of scales enters even for ideals. This is described in section 4.2 of \cite{G5}.

The effect of considering the multiplicity of the pair is to return to the situation of modules of finite colength, where there is only a single stratum, which is of dimension $0$. Because of this, in the sequel to this paper we are able to prove a principle of specialization of integral dependence, which generalizes the main result of \cite{GK}.

In section 1, we give the main result of the paper. Besides the notion of the multiplicity of a pair of modules, a notion due to S. Kleiman and A. Thorup, we use the notion of a mixed polar of a pair of modules $(M,N)$. In section 2, we describe some applications of the main result of section 1. 
The applications of the multiplicity-polar theorem are of two types. In the first we construct a deformation so that we can understand the significance of the multiplicity of the pair, or show we can use the multiplicity of the pair to describe some geometric behavior. Applications of this type in section 2 are the geometric significance of the BR-mulitplicity (Theorem 2.1), the geometric significance of multiplicity of the pair 
$(J(f), I)$ where $f$ defines a hypersurface with non-isolated singularities, $J(f)$ is the jacobian ideal of $f$ and $I$ is the ideal of the singular locus of the hypersurface. We do the cases where $I$ defines a complete intersection (Theorem 2.3), or the hypersurface is the image of a finitely determined map-germ $F:\IC^2,0\to \IC^3,0$ (Theorem 2.6). We then show how the multiplicity of the pair can be used to calculate the index of a differential 1-form with an isolated singularity on a germ of an isolated singularity (Theorem 2.8). In the second kind of application, we use the generalization of the Principle of Specialization of Integral Dependence to prove results in equisingularity. The particular condition we study in this paper is the \Wf \hskip 2pt condition (Theorems 2.10 and 2.11.)

The author thanks Steven Kleiman for many helpful conversations during the period of years in which this work developed.

\sctn  The Multiplicity-Polar Theorem

Given a submodule $M$ of a free $\O_X^d$ module $F$ of rank $p$, we can associate a subalgebra ${\cal R }(M)$ of the symmetric
$\O_X$ algebra on $p$ generators. This is known as the Rees algebra of $M$. If $(m_1,\dots,m_p)$ is an element of $M$
 then $\Sum m_iT_i$ is the corresponding element of 
${\cal R}(M)$. Then $\Projan ({\cal R} (M))$, the projective analytic spectrum of ${\cal R}(M)$ is the closure of 
the projectivised row spaces of $M$ at points where the rank
of a matrix of generators of $M$ is maximal. Denote the projection to $X$ by $c$, or by $c_M$ where there is ambiguity. We call $\Projan ({\cal R} (M))$ {\it the $M$-transform} of $X$. In the development of the theory of the integral closure of modules it sometimes plays a role similar to that of the blowup of a space by an ideal. Let $C(M)$ denote the locus of points where $M$ is not free, ie. the points where the rank of $M$ is less 
than $e$; call it the critical locus of $M$, let $C(P(M))$ denote $c^{-1}(C(M))$.

Throughout this paper we assume $X^d$ is an equidimensional,  analytic space, reduced off a nowhere dense subset  and,  we assume the generic rank of $M$ is $e$ on each 
component of $X$. The hypothesis on the equidimensionality of $X$ and on the rank of $M$ ensures that $\Projan{\cal R}(M)$ 
is equidimensional of 
dimension $d+e-1$.

If $M$ is an ideal, then the $M$-transform of $X$ is just the blowup of $X$ by $M$, and the $M$ critical locus is just the zero set of $M$. Denote the $M$-transform of $X$ by $T_M(X)$, $c_M(X)$, the $M$ critical set.

Suppose $(X,x)$ is a complex
analytic germ, $M$ a submodule of
$\O_{X,x}^{p}$. Then $h\in \O_{X,x}^{p}$ is in the {\it{integral closure of
$M$}}, denoted $\?{M}$,
if and only if for all
$\phi :\left({\IC},0\right) \rightarrow \left(X,x\right)$, $ h \circ \phi
\in (\phi^{*}M){\O}_{1}$. (For the development of the theory of the integral closure of modules from this point of view and its connection with the development by Rees \cite {R}, see \cite{G2}.)

If $N$ is a submodule of $M$ and $\?M=\?N$ we say that $N$ is a {\it reduction} of $M$.

In this paper we will want to use the multiplicity
of a pair of modules $M\subset N$, $\?M$ of finite colength in $\?N$, as well. 
We recall how to construct these numbers following the approach
 of Kleiman and Thorup (\cite{KT1}).
 
  If $M$ is a submodule of $N$, then ${\cal R}(M)$ is a subalgebra of ${\cal R}(N)$. Denote the ideal generated by the elements of  $\cal R(M)$  of degree 1 by $\cal M$. If $M$ is a submodule of $N$or $h$ is a section of $N$, then
 $h$ and  $M$  generate ideals on $\Projan{\cal R}(N)$; denote them
by $\rho(h)$ and $\rho({\cal M})$. If we can express $h$ in terms of a set of generators $\{n_i\}$ of $N$ as $\Sum g_in_i$, 
then in the chart in which $T_1\ne 0$, we can express a generator of $\rho(h)$ by $\Sum g_iT_i/T_1$. 
The effect of passing to the ideal sheaf $\rho({\cal M})$ is to divide $M$ by $N$.  In order to study the remainder, which is where the difference between $M$ and $N$ lies,  we blowup the $N$-transform of $X$ by 
$\rho({\cal M})$.

On the blowup $B_{\rho({\cal M})}(\Projan {\cal R}(N))$ we have two tautological bundles, one the pullback of 
the bundle on $\Projan{\cal R}(N)$,
 the other coming from $\Projan{\cal R}(M)$; denote the corresponding Chern classes by $l_M$ and $l_N$. 
Suppose the generic 
rank of $N$ (and hence of $M$) is $e$.
 Then 
the multiplicity of a pair of modules $M,N$, as developed by Kleiman and Thorup (\cite {KT1}, is:

$$e(M,N)=\Sum_{j=0}^{d+e-2} \int D_{M,N} \cdot l^{d+e-2-j}_{M}\cdot l^j_{N}.$$

If $M\subset \O^p_X$ is of finite colength, then the multiplicity of $M$ is the multiplicity of the pair $(M,\O^p_X)$.

At this point we discuss how the generalization of the Buchsbaum-Rim multiplicity along the lines of \cite{G-G} would go. Notice that if $M$ is a submodule of a free module $F$ of rank $p$, and  the generic rank of $M$ is $p$, then the above intersection number is well defined with $N$=$F$, even though $M$ does not have finite colength in $F$. This number can be taken as the generalization to modules of the zero dimensional Segre number of ideals defined in  \cite{G-G}. Given a prime ideal $P$ in $\O_X$, we can choose a generic plane slice $H$  at a generic point of $V(P)$ of complementary dimension to $V(P)$, restrict $M$, $F$ to this slice, then calculate the analogous intersection number. This  gives a function on prime ideals. It is not hard to see that this number is zero unless $V(P)$ is the image of a component of the exceptional divisor of the blowup by $\rho({\cal M})$. At this point there are two directions to go. In one, you multiply the value of the function on a prime by the multiplicity of the prime and sum over all primes of the same height. This gives a single sequence of numbers which specialize in the ideal case to the Segre numbers of  \cite {G-G}.  The other direction is to simply stay with the function on the primes.

In order to study how $e(M,N)$ changes in a family and to describe the problem associated with the above notions of multiplicity, we need the notion of the polar varieties of $M$ and $N$ and also that of their mixed polar varieties.

The polar variety of codimension $k$ of $M$ in $X$ denoted $\Gamma_k(M)$ is constructed by intersecting $\Projan{\cal R}(M)$
 with $X\times H_{e+k-1}$ where 
$H_{e+k-1}$ is a general plane of codimension $e+k-1$, then projecting to $X$. 

This notion was developed by Teissier 
in the case where 
$M=JM(F)$, the jacobian module of $X$.  (\cite{T-2}). (If $X^d=F^{-1}(0$), then the jacobian module of $X$ is the module generated by the partial derivatives of $F$).    Here, think of $H$ as the projectivised row space of a linear submersion $\pi$ mapping $\IC^n\to \IC^{d-k+1}$, where $X^d\subset\IC^n$.  Note that $\Projan{\cal R}(JM(F))\subset X\times \IP^{n-1}$. The dimension of $H$ is $d-k$, so its codimension is $(n-1)-(d-k)=(n-d)+k-1$ as it should be. The intersection of the $JM(F)$-transform of $X$ with  $X\times H_{e+k-1}$ generically consists of pairs $(x,P)$ where $x$ is a smooth point of $X$ and $P$ a hyperplane which contains the kernel of  $D(F)(x)$ and $\pi$.
 Thus, a Zariski open and dense subset of $\Gamma_k(JM(F))$ consists of
 the set 
of smooth points of $X$ where the matrix formed from $D\pi$ and $DF$ has less than maximal rank possible, which is $e+d-k+1=n-k+1$ and hence greater than minimal kernel
 rank. These
are just the points where
the restriction of $\pi$ to the smooth part of $X$ is singular. 

In dealing with $\Projan{\cal R}(M)$ it is helpful to denote the number of generators of M by $p_M$.

In general, we can think of $\Gamma_k(M)$ as the closure of the set of points where the module $M$ has maximum rank and 
where the module whose matrix of generators consists of the matrix of generators of $M$ augmented by the rows of 
the linear submersion
$\pi:\IC^{p_M}\to \IC^{p_M-k-e+1}$, has less than maximal rank $p_M-k+1$. 

When we consider $M$ as part of a pair of modules $M,N$, where the generic rank of $M$ is the same as the generic rank of $N$,
 then other polar varieties become interesting as well. In brief, we can intersect
 $B_{\rho({\cal M})}(T_N(X))\subset X\times \IP^{p_M-1}\times\IP^{p_N-1}$ with a mixture of hyperplanes from the two projective spaces
which are factors of the space in which the blowup is embedded. We can then push these intersections down to 
$T_N(X)$ or
$X$ as is convenient, getting mixed polar varieties in $T_N(X)$ or in $X$. Denote by $\Gamma_{i,j}(M,N)$ the mixed polar variety on $X$ defined using $i$ hyperplanes from  $\IP^{p_M-1}$ and $j$ hyperplanes from $\IP^{p_N-1}$. Denote the corresponding mixed polar on $T_N(X)$ by  $\Gamma_{i,j}(M,N)_{T_N}$. Notice that 
$\Gamma_{i,0}(M,N)=\Gamma_{i-(e-1)}(M)$ for $i>(e-1)$, and in general the codimension of $\Gamma_{i,j}(M,N)$ is $i+j-(e-1)$ for $i+j>(e-1)$.

Now we extend these ideas to families of sets and modules.

Setup: We suppose we have families  of modules $M\subset  N$, $M$ and $N$
 submodules of a free module $F$ of rank $p$
 on an equidimensional family of spaces with equidimensional
 fibers ${\cal X}^{d+k}$, ${\cal X}$ a family over a smooth base
$Y^k,0$. We denote the fiber of the family over $y\in Y$ by $N(y)$, $M(y)$, $X(y)$. We assume that the generic rank of $M(y)$, $N(y)$ is $e\le p$ on each component of $X(y)$ for all $y$ in  neighborhood of $0$.  
Let $C({\cal M})$ denote the cosupport of $\rho({\cal M})$ in $T_N{\cal X}$. 

Notice that the various intersection numbers which appear in the formula for the multiplicity of the pair can be calculated by finding the degree of $\rho(M)$ restricted to the curves $\Gamma_{i,j}(M,N)_{T_N}$, where $i+j=d+e-2$. In order to see how the multiplicity of the pair changes in the family, we will see how these intersection numbers change. To do this we need a lemma which 
shows that the various mixed polars  specialize properly.

\lem 1 Suppose in the above setup we have that $\?M=\?N$ off a set of dimension $k$ which is finite over $Y$. Suppose further
that $C(P(M))(0)=C(P(M(0)))$ except possiby at the points which project to $0\in {\cal X}(0)$. Then all of the mixed polars of $M$ and $N$ of dimension
$k+1$ specialize over $0\in Y$.

\pf The hypothesis  $C(P(M))(0)=C(P(M(0)))$ except possiby at the points which project to $0\in {\cal X}(0)$
 means that if a component of 
$C(P(M))(0)$ projects to a set of dimension $j>0$, then the  dimension of the  fiber of the component 
over $0\in {\cal X}(0)$ is at most $d+e-3$, since the dimension of any component
of $C(P(M(0)))$, hence of $C(P(M))(0)$ is at most $d+e-2$, and the image of the components we consider
have dimension at least $1$.

This means that if we intersect $T_M({\cal X})$ with $d+e-2$ generic hyperplanes, the intersection of the hyperplanes
 will miss the fiber of these components over $0\in X(0)$, hence in a small enough neighborhood
of $C(P(M))(0,0)$, will miss the components.
 This implies
that the fiber of the polar of $M$  of dimension
$k+1$ specializes over $0\in Y$. (For the failure to specialize implies that some component of the fiber over $0\in Y$ lies in 
a component of $C(P(M))(0)$ which projects to a set of positive dimension in $X(0)$. That this component is unavoidable
implies that its fiber over
$0\in {\cal X}(0)$ has dimension at least $d+e-2$.)

Since $\?M=\?N$ off a set $K$ of dimension $k$ which is finite over $Y$,  $T_N({\cal X}-K)$ is finite over $T_M({\cal X}-K)$.
This implies that the  dimension of the  fiber of any component of $C(P(N))(0)$
over $0\in {\cal X}(0)$ which projects to a set of dimension $j>0$, is at most $d+e-3$ also.  
So the polar of $N$ specializes as well.

Now for the mixed polars.

Recall that we have a birational map from $T_N{\cal X}$ to $T_M{\cal X}$ because $M$ and $N$ have the same generic rank, hence 
both spaces are birational with
$B_{\rho({\cal M}))}(T_N{\cal X})$ (\cite{KT1}), and this map is an equivalence off $C({\cal M})=V(\rho(M))$. This means that the polar of 
$M$ of dimension $k+1$ can also be constructed by setting $d+e-2$ elements of $\rho({\cal M})$ equal to zero in $P(N)$ 
slashing out by the cosupport of 
$\rho({\cal M})$ and projecting to ${\cal X}$. To construct the mixed polar with $i$ elements from $M$ and $j$ elements from $N$, 
intersect $T_N{\cal X}$ with $j$ generic hyperplanes, restrict $\rho({\cal M})$ to this intersection, and form the polar of codimension $i$,
 then project to 
${\cal X}$.  Suppose $\{m_i\}$, $1\le i\le d+e-2$ define the polar of dimension $k+1$ of $M$. 
If we drop $m_{d+e-2}$ from this collection, 
assume the remainder of the collection defines the polar of dimension $k+2$. Then the components of the zeroes of the remainder
which do not lie in $V(\rho({\cal M}))$ can intersect $T_N{\cal X}(0)$  in at most a surface, and the surface intersects $V(\rho({\cal M}))$ in at most a curve,
and the fiber of the components of $C(P(N))(0)$ which project to sets of positive dimension  in $X_0$ in a  finite number
of points.   Then we can choose a hyperplane which misses these points and curves, so the 
mixed polar defined using this hyperplane
and $\{m_i\}$, $1\le i\le d+e-3$ specializes properly.  A similar argument works for
 the other polars as well. This finishes the proof.

\lem  2  The projections to ${\cal X}$ and  to ${\cal X}(y)$
induce birational maps between 
 $\Gamma_{i,j}(M,N)_{T_N}$ and $\Gamma_{i,j}(M,N)$ and between $\Gamma_{i,j}(M(y),N(y))_{T_N(y)}$ and  $\Gamma_{i,j}(M(y),N(y))$)  for $i+j\ge e-1$.
 
 \pf It is clear that  the image of the projection is $\Gamma_{i,j}(M,N)$. At a general point $x$
 of this set $M$ and $N$ have rank
$e$, so the dimension of the fiber of $T_N{\cal X}\cap H_j$ over $x$ is $e-1-j$, $H_j$ a codimension $j$ plane
on $T_N{\cal X}$. The zeroes of the $i$ generic elements of $\rho(M)$  on the 
fiber of  $T_N{\cal X}\cap H_j$ over $x$ are another linear subspace of codimension $i$, so in general
 the intersection of these two spaces and the fiber of $P(N)$ over $x$ is a point, since the sum of their codimensions is $\ge e-1$.

 Now we are almost ready to prove a first result relating intersection numbers of pairs and mixed polars. Given a mixed polar on ${\cal X}$ of dimension $k$ passing through the origin, we can consider the degree of the map from the polar to our smooth base $Y$ at the origin. We denote this number by 
 ${\rm mult}_Y \Gamma_{d+e-1-j,j}(M,N)$ where $\Gamma_{d+e-1-j,j}(M,N)$ is the polar in question. In the following Theorem, part of being a generic point of $Y$ is that the mixed polars of dimension $k$ are empty at $y$.
 
 \thm 3 Suppose in the above setup we have that $\?M=\?N$ off a set $C$ of dimension $k$ which is finite over $Y$. 
Suppose further
that $C(P(M))(0)=C(P(M(0)))$ except possiby at the points which project to $0\in {\cal X}(0)$. 
Then, for $0\le j\le e+d-2$ and $y$ a 
generic point of $Y$,

$$\int D_{M(0),N(0)} \cdot l^{d+e-2-j}_{M(0)}\cdot l^j_{N(0)}
-\int D_{M(y),N(y)} \cdot l^{d+e-2-j}_{M(y)}\cdot l^{j}_{N(y)}=$$

$${\rm mult}_y \Gamma_{d+e-1-j,j}(M,N)-{\rm mult}_y \Gamma_{d+e-2-j,j+1}(M,N).$$

\pf We work on $T_N({\cal X})$, which contains  $T_N({\cal X}(0))$and  $T_N({\cal X}(y))$. 

In Lemma 1.1 we showed that we could pick elements defining the mixed polars $\Gamma_{i,j}(M,N)_{T_N}$, where $i+j=d+e-2$, such that $\Gamma_{i,j}(M,N)_{T_N}(0)=\Gamma_{i,j}(M(0),N(0))_{T_N(0)}$.  The next step is to show, that at points of $C(y)$ for generic $y$ that $\Gamma_{i,j}(M,N)_{T_N}(y)=\Gamma_{i,j}(M(y),N(y))_{T_N(y)}$.

Note that this is not true if $C$ is not finite over $Y$ even if $M$ is an ideal and $N=\O_{\cal X}$. For example if 
dim ${\cal X}$ is 3, dim $Y$ is 1, then the exceptional divisor 
of $B_I({\cal X})$ may have a component $V$ with generic fiber over its image
 of dimension 0 and dim $V(0)=1$, but with fiber over the origin of dimension 1. Then a generic hyperplane must intersect
the fiber over zero, hence intersect the component in a curve, which we can take to be finite over $Y$. This implies that $\Gamma_1(I)$ which is a surface will 
contain curves in the cosupport of $I$ such that the polar curve of $I_y$ at any $x$ on the curve $\Gamma_1(I)(y) $is empty, 
in particular the fiber of $\Gamma_1(I)$ over $y$ is not a polar curve for $I_y$ at $x\in C(y)$.

We begin with $\Gamma_{d+e-2,0}(M,N)$. We know that $ V(\rho(M))$ is a collection of varieties which project to varieties of dimension
at most $k$. Consider $B_{\rho(M)}(T_N{\cal X})$. 
If we work over a generic point $y\in Y$, the components of its exceptional divisor $D$ which have 
non-empty fiber over $y$ must surject onto $Y$, hence onto a component of $C$, since 
these are finite over $Y$.
 By \cite{KT2} p402 the image of each of the components of $D$ in $T_M({\cal X})$ has 
dimension $e+d+k-2$, and 
each such image projects to a component of $C\subset{\cal X}$ of dim $k$,
 so the fiber dimension is at least $d+e-2$. 
Now generically $T_M{\cal X}(y)=T_{M(y)}{\cal X}(y)$, so the fiber dimension of $T_{M(y)}{\cal X}(y)$
over each point of $C(y)$ is $d+e-2$.
This implies that the polar curve of $M(y)$ at such a point is non-void. Further 
by varying our elements from ${\cal M}$ and hence our hyperplanes on $T_M{\cal X}$, we can ensure that our
 hyperplanes intersect the image of
each component of $D$ transversely generically. Picking a $y$ over which our setup satisifies all of the genericity conditions, 
if necessary moving our hyperplanes slightly, we can ensure that the germ of $\Gamma_{d+e-2,0}(M(y),N(y))$ is a  
polar curve at each of the points in the fiber over $y$ where the integral closure of
$M(y)$ is different from $N(y)$. (Slight moves of our hyperplanes do not move the polars off 
these points from the generic transversality conditions.)

For $\Gamma_{d+e-2-j,j}(M,N)$, we choose j hyperplanes on $T_N{\cal X}$ for the elements which 
contribute the j to the mixed polar,
 whose intersection is transverse to the fibers of $V({\cal M})$ over $0$, and such that on 
$B_{\rho(M)}(T_N{\cal X})$ the image 
of its intersection with the components of  $D$ over $0$ has dimension $d+e-2-j$. 
Now each component
of $D_j$ the exceptional divisor of $B_{\rho(M)}(T_N{\cal X}\cap H_j)$ will map to $T_M{\cal X}$  with 
dimension $d+e+k-2-j$, hence,
given any set of hyperplanes on $T_M{\cal X}$ with $d+e-2-j$ elements, 
the intersection of the hyperplanes will intersect 
each of the fibers over $Y$
of the images of the components of $D_j$.  Note further, that
 only components of $V({\cal M})$ of dimension $\ge k+j$ will 
give a component of $V({\cal M})\cap H_j$ which surjects onto $Y$.
 Since we are working over a generic point of $Y$,
 we can ignore all other components.

Suppose $x\in {\cal X}$ is in the image of $V({\cal M}(y))\cap H_j$, $x$ over a generic $y\in Y$.
 Then if we move $H_j$,
 each component of $V({\cal M}(y))$ which lies over $x$ and meets $ H_j$ will continue to meet $H_j$
because of the dimension condition on these components. By the argument of the previous paragraph if we move $H_j$, and our 
$d+e-2-j$ elements of ${\cal M}$, the mixed polar $\Gamma_{d+e-2-j,j}(M,N)_{T_N}$ will continue to meet each component
of $(V{\cal M}(y))$ which lies over $x$.  Again by slight moves of our hyperplanes we can ensure that  
$\Gamma_{d+e-2-j,j}(M,N)(y)$ is  $\Gamma_{d+e-2-j,j}(M(y),N(y))$ at $x$.

Now we choose $R\in{\cal M}$ and $Q$ a linear form on $T_N{\cal X}$, which satisfy the following conditions. If we add $R$ to the elements defining $\Gamma_{d+e-2-j,j}(M,N)_{T_N}$ we get 
$\Gamma_{d+e-1-j,j}(M,N)_{T_N}$. Since the fiber of  $\Gamma_{d+e-2-j,j}(M(y),N(y))_{T_N(y)}\cap V({\cal M})$ 
is also zero dimensional, we can choose $Q$ so that the zeroes of $Q$ 
generically miss these points and  if we add $Q$ to the set of elements defining $\Gamma_{d+e-2-j,j}(M,N)_{T_N}$, we get $\Gamma_{d+e-2-j,j+1}(M,N)_{T_N}$. Further, $R/Q\in \rho(M)$ is a reduction of 
$\rho(M)(0)$ on $\Gamma_{d+e-2-j,j}(M(0),N(0))_{T_N(0)}$, 
and is a reduction of $\rho(M)$ restricted to the polar $\Gamma_{d+e-2-j,j}(M,N)_{T_N}$ on a non-empty
Z-open dense subset $V({\cal M})$. (This set may not include the fiber 
of $\Gamma_{d+e-2-j,j}(M,N)_{T_N}$ over $0$, regarded as a subset of $\Gamma_{d+e-2-j,j}(M,N)_{T_N}$.
$R/Q$ may only be a reduction of $\rho(M)$ restricted to 
 $\Gamma_{d+e-2-j,j}(M(0),N(0)_{T_N(0))}$ .) This is true because the fibers of 
$\Gamma_{d+e-2-j,j}(M,N)_{T_N}$  form a family of curves over $Y$; so the element $R/Q$ gives the
desired reduction at points of $V({\cal M})$ which do not lie on the polars of dimension $k$ defined by $R$ and $Q$.

Then we have at a generic $y$ value:

$$\int D_{M(y),N(y)} \cdot l^{d+e-2-j}_{M(y)}\cdot l^{j}_{N(y)}$$
$$=\sum{\rm deg} R/Q|(\Gamma_{d+e-2-j,j}(M(y),N(y))_{T_N(y)},x),$$
where the second sum is taken over points $x\in C({\cal M}(y))$. This equality also holds at $y=0$, where all of the points $x$ are in the fiber of $C({\cal M})(0)$ over the origin in $X(0)$. Since there are only a finite number of such $x$ in the polar curve $\Gamma_{d+e-2-j,j}(M(0),N(0))_{T_N(0)}$, the sum is well defined.

Now we are ready to work on ${\cal X}$.

Because $\pi_{\cal X}$ restricted to $\Gamma_{d+e-2-j,j}(M,N,P(N))$ is birational it follows that 
there are functions $R'$ and $Q'$
on $\Gamma_{d+e-2-j,j}(M,N)$ such that $R'/Q'\circ \pi_{\cal X}=R/Q.$ Let $R_0$, $R_y$, $Q_0$, $Q_y$ denote the restrictions to the fibers over $0$, $y$.

Now we have:

$$\int D_{M(0),N(0)} \cdot l^{d+e-2-j}_{M(0)}\cdot l^j_{N(0)}
-\int D_{M(y),N(y)} \cdot l^{d+e-2-j}_{M(y)}\cdot l^{j}_{N(y)}=$$

$$\hskip -1in{\rm deg} R'_0/Q'_0|\Gamma_{d+e-2-j,j}(M(0),N(0))-$$
$$\Sum_{x\in{\pi^{-1}_Y(y)\cap C}}
 {\rm deg} R'_y/Q'_y|(\Gamma_{d+e-2-j,j}(M(y),N(y)),x).$$

More is true.

The degrees of $(\pi_y,R')$, $(\pi_y,Q')$ are locally constant; 
since the fiber over zero of $\pi_Y$ is generically reduced,
it follows that the degree of $\pi_y,R'$ at $0\in {\cal X}$ equals the degree of 
$R'_0$ restricted to $\Gamma_{i,j}(M(0),N(0))$ and a similar statement holds for $Q'$. So the number of zeroes of $R'_y/Q'_y$ less the number of poles of 
 $R'_y/Q'_y$ is equal to a constant. Denote the  zeroes of $R'_y$ by $Z(R'_y)$.  Putting together the above ideas we get:
 
 $$\hskip -1in{\rm deg} R'_0/Q'_0|\Gamma_{d+e-2-j,j}(M(0),N(0))=$$
 $$\#(Z(R'_y)/C(y))-\#(Z(Q'_y)/C(y)= $$
 $$\Sum_{x\in{\pi^{-1}_Y(y)\cap C}}
 {\rm deg} R'_y/Q'_y|(\Gamma_{d+e-2-j,j}(M(y),N(y)),x)+{\rm  mult}_y \Gamma_{d+e-1-j,j}(M,N)$$
 $$-{\rm mult}_y \Gamma_{d+e-2-j,j+1}(M,N).$$
 
 Hence

$$\int D_{M(0),N(0)} \cdot l^{d+e-2-j}_{M(0)}\cdot l^j_{N(0)}
-\int D_{M(y),N(y)} \cdot l^{d+e-2-j}_{M(y)}\cdot l^{j}_{N(y)}=$$
$${\rm  mult}_y \Gamma_{d+e-1-j,j}(M,N)-{\rm mult}_y \Gamma_{d+e-2-j,j+1}(M,N).$$

\cor 4 (Multiplicity-Polar Theorem) With the hypotheses of Theorem 1.3
$$\Delta e(M,N)={\rm  mult}_y \Gamma_{d+e-1,0}(M,N)-{\rm mult}_y \Gamma_{0,d+e-1}(M,N)$$

$$={\rm  mult}_y \Gamma_{d}(M)-{\rm mult}_y \Gamma_{d}(N).$$

\pf By \cite{KT1} p191,
$$e(M(0)),N(0))=\Sum_{j=0}^{d+e-2} \int D_{M(0),N(0)} \cdot l^{d+e-2-j}_{M(0)}\cdot l^j_{N(0)}$$

So, summing the equalities of 1.3 over $j$, and using the relations between the mixed polars of type $(i,0)$, $(0,j)$ and the polar varieties of $M$ and $N$ we obtain the desired result.

Although we postpone the most general version of the principle of specialization of integral dependence until the sequel to this paper, we prove an elementary version in the setting of 1.3, to show how  Theorem 1.3 comes into play in this question.

First we relate $T_M({\cal X})$ and Corollary 1.4.

\cor 5  Suppose with the hypotheses of Theorem 1.3, that either 

1) $\Delta e(M,N)=0$ and $\Gamma_{d}(N)$ is empty, or 

2) $e(M(0),N(0) +{\rm mult}_y \Gamma_{d}(N)=e(M(y),N(y))$, $y$ a generic point of $Y$.

Then 

1) $\Gamma_{d}(M)$ is empty.

2) The dimension of the fiber of $T_M({\cal X})$ over the origin in ${\cal X}$  is less than $d+e-1$.

\pf  By Corollary 1.4 either 1) or 2) imply that  $\Gamma_{d}(M)$ is empty, since it must have multiplicity $0$. If  $\Gamma_{d}(M)$ is empty, then the intersection of $T_M({\cal X})$
with $H_{d+e-1}$ must be empty. Since $H$ is general, this intersection is empty if and only if the intersection in $\IP^{g(M)-1}$ of $H$ and  the fiber of $T_M({\cal X})$ over the origin is empty. So the dimension of the fiber of $T_M({\cal X})$ over the origin is less than $d+e-1$.

\cor 6 With the hypotheses of Theorem 1.3, suppose $h\in F$, and that there exists a Zariski open subset $U$ of $Y$ such that $h$ is in $\?{M_x}$ for all $x\in {\cal X}|U$. Then if either of the hypotheses of Corollary 1.5 hold, $h\in\?M$.

\pf 
The main tool we use is a result of Kleiman and Thorup which relates the dimension of the fiber of $T_M(X)$ and the integral closure of $M$. The result is 
Theorem A.1 of \cite{KT2}. It says that given two modules $M\subset N\subset  F$, $F$ free, $N$ and $M$ generically equal and free of rank $e$, and $W$ the set of points on $X^d$ where $N$ is not integral over $M$, then the dimension of $T_M(X)|_W$ must be $d+e-2$.


 Consider the module $N=M+(g)$. Let $W$ be the closed subset where $N$ is not integral over $M$, it is clear that $W\subset C(M)$.
 Set 
$E:=c^{-1}(W)$, $c$ the strucure map of $T_M(X)$, so $E\subset C(P(M))$. If $V$ is a component of  $C(P(M))$, there are two cases--$c(V(0))$ is the origin or it is not. If it is not, then by the hypothesis of 1.3, $V(0)\subset C(P(M(0)))$, hence has dimension less than $d+e-1$. Since the fiber dimension of $V$ over the origin is less than $d+e-1$, the dimension of $V|_W$ is less than $(d+e-1)+(k-1)$. This implies that any component of $E$ in $V$ also has dimension less than $d+k+e-2$. Suppose $c(V(0))=0$. Then by 1.5, the dimension of $V(0)$ is also less than $d+e-1$, so any component of $E$ in such a $V$ again has dimension less than $d+k+e-2$. Then the result of Kleiman and Thorup implies that $W$ is in fact empty, so $g\in \?M$, which finishes the proof.

Earlier in this section we described how the generalization of the Buchsbaum-Rim multiplicity along the lines of \cite{G-G} would go. As  \cite{G6} shows, the invariants in this approach in a 1-dimensional family may be independent of parameter, yet there may still be a polar curve, hence the dimension of 
the fiber of $T_M({\cal X})$ over the origin in ${\cal X}$  is  $d+e-1$, which means that the theorem of Kleiman and Thorup does not apply.

\sctn  Applications

The applications of the multiplicity-polar theorem are of two types. In the first we construct a deformation so that we can understand the significance of the multiplicity of the pair, or show we can use the multiplicity of the pair to describe some geometric behavior.

The second type of application is based on the ability of the multiplicity-polar theorem to relate infinitesimal information, such as the change in a family of the multiplicity of a pair of modules, and macroscopic information such as the polar varieties of the modules. Stratification conditions such as the Whitney conditions or Thom's \Af condition can be phrased in terms 
of limits of linear spaces. These conditions in turn can be controlled by the theory 
of integral closure of modules. By using the multiplicities of pairs of modules it is possible to verify these conditions stratum by stratum; the multiplicity-polar theorem is the key tool in showing that the independence of our invariants from parameter implies that the polar curve of our module is empty, which in turn implies the integral closure condition we need.

 We first describe an application of the first type  to the simple case where $N$ 
is free.

As mentioned in the introduction, the following geometric interpretation of the multiplicity of an ideal is 
well known, and lies at the heart of the application of multiplicity to intersection theory:
Given an ideal $I$ of finite colength in $ \O_{X,x}$, $X^d$ 
equidimensional, choose $d$ elements $(f_1,\dots,f_d)$ of 
$I$ which generate a reduction of $I$. 
 Then the multiplicity of $I$ is the degree at $x$ of $F$ where $F$
is the branched cover defined by the map-germ with components  
$(f_1,\dots,f_d)$, 
or the number of points in a fiber of $F$ over a regular value close to 
$0$. (Cf. \cite{M} p )

We wish to give a similar interpretation of the multiplicity of a module. 

\thm 1 Given $M$ a submodule of  $\O^p_{X,x}$, $X^d$ equidimensional, choose
$d+p-1$ elements which generate a reduction $K$ of $M$. Denote the matrix 
whose columns are the $d+p-1$ elements by $[K]$; $[K]$ induces
a section of $\Hom \hskip 2pt(\IC^{d+p-1},\IC^p)$ which is a trivial 
bundle over $X$. Stratify $\Hom \hskip 2pt(\IC^{d+p-1},\IC^p)$ by rank. 
Let $[\epsilon]$ denote a $p\times (d+p-1)$ matrix, whose entries are 
small, generic constants. Then, on a suitable neighborhood $U$ of $x$ 
the section of 
$\Hom \hskip 2pt(\IC^{d+p-1},\IC^p)$ induced from $[K]+[\epsilon]$ has at 
most kernel rank 1, is transverse to the rank stratification, 
and the number of points where the kernel rank is 1 is $e(M)$.

\pf The first step is to explain by construction what we mean by ``generic 
constants". Consider the family of maps $G_a$ from 
$X^d$, parametrised by $\IC^{p(d+p-1)}$ to
$\Hom \hskip 2pt(\IC^{d+p-1},\IC^p)$ defined by 
$G_a(x)=G(x,a)=[K(x)]+[A]$, where $[A]$ is the $p\times (d+p-1)$ matrix 
whose entries are coordinates
$a_{i,j}$ on $\IC^{p(d+p-1)}$. Let $\tilde X$ be a resolution of $X$, so 
we have an induced family of maps $\tilde G$ on $\tilde X$. 
Since the map $\tilde G(x,a)$ is a submersion, it follows that for a 
Z-open subset $V$ of $\IC^{p(d+p-1)}$, that for $a\in V$, the map
$\tilde G_a$ is transverse to the rank stratification. We claim that the 
points of $V$ are the generic constants in the theorem.
 Note that the points of $\Hom \hskip 2pt(\IC^{d+p-1},\IC^p)$ of kernel 
rank 1 have codimension $1\cdot((d+p-1)-(p-1))=d$; so since 
$\tilde G_a$ is transverse it can only hit points of the rank 
stratification of kernel rank 1, and only if $D\tilde G_a$ has maximal 
rank at such points which implies
$X$ is smooth at the projection of such points. Let $\tilde K$ be the 
submodule of ${\cal O}^p_{X\times \IC^{p(d+p-1)}}$ defined by the matrix
$[K(x)]+[A]$.

Now we apply the multiplicity-polar theorem to $X\times \IC^{p(d+p-1)}$, 
thought of as a family parametrised by $\IC^{p(d+p-1)}$, and 
$(\tilde K,\O^p_{X\times \IC^{p(d+p-1)}})$. Use a point of $V$ as the 
generic parameter value $\epsilon$. 
Then $\O^p_{X\times \IC^{p(d+p-1)}}$ has no polar, because it is free, 
$\tilde K$ has no polar,
 because
$\tilde K$ is generated by $d+p-1$ elements.  (Hence if all these elements vanish we are already on the cosupport of $\tilde K$.)

The hypothesis that $C(P(\tilde K))(0)=C(P(\tilde K(0)))$  except possibly over $x$, holds because $\tilde K(0)=M$ is free in a neighborhood of $x$ in $X$ except possibly at $x$.

Choose $U$ a neighborhood of 
$x\times \IC^{p(d+p-1)}$ sufficently small such that every component
 of the cosupport of $\tilde K$ which meets $U$ has $(x,0)$ in its 
closure. Now at $\epsilon$ the cosupport of $\tilde K_{\epsilon}$ is just 
the points where
  $[K]+[\epsilon]$ has less than maximal rank. At such points $e(\tilde 
K_{\epsilon})$ is $1$, because since we are at a smooth point of $X$, the 
local ring of $X$ is Cohen-Macaulay,
so $e(\tilde K_{\epsilon})$ is just the colength, which is 1. Hence 
$e(M)=e(K)=e(\tilde K_{0})=e(\tilde K_{\epsilon})$, which is the number of 
points 
 where the kernel rank of $[K]+[\epsilon]$ is 1.

We now return to the setting of the introduction.

 \cor 2 Suppose $X^d$ is an equidimensional space, $E$ a vector bundle on $X$ of rank $p$, $M$ a submodule of the sheaf of sections of $E$, and in a neighborhood of a point $x\in X$ the rank of $M$ is $p$. Suppose $N$ is a subsheaf of $M$ generated by $d+p-1$ sections such that $N$ is a reduction of $M$ at $x$. Then the number of points where a generic perturbation of the generators of $N$ has less than maximal rank is $e(M,x)$.

\pf Locally, the sheaf of sections of $E$ is a free module on $p$ generators. The result then follows from Theorem 2.1.

 In \cite{F} p254 Fulton describes the k-th degeneracy class 
associated to $\sigma$ a homomorphism of vector bundles over $X^d$. The 
support of the class is the set of points where
the rank of $\sigma$ is less than or equal to $k$. Suppose $\sigma: E\to 
F$ where the rank of $E$ is $e$ and the rank of $F$ is $f$, $e\ge f$, 
$e-f+1=d$. Then the $f-1$ degeneracy class is 
supported at isolated points. Fulton shows that if $X$ is Cohen-Macaulay 
at $x$, the contribution to the class at $x$ is the colength of the ideal 
of maximal minors of the
matrix of $\sigma$ at $x$  for some suitable local trivializations of $E$ 
and $F$. Note that this is just the Buchsbaum-Rim multiplicity of the 
module generated by the columns
of the matrix  associated to $\sigma$. Theorem 1.2 shows that in this 
situation if $X$ is pure dimensional, 
the contribution to the degeneracy locus is always the Buchsbaum-Rim 
multiplicity associated to $\sigma$ at $x$, so the Cohen-Macaulay hypothesis 
is unnecessary. (Just use the proof of 1.2 to construct a rational 
equivalence to go back to Fulton's case 
close to $x$.) 

The next application occurs in the study of families of hypersurfaces with $1$-dimensional singular locus. Here is the setup for families of functions defining hypersurfaces. Assume that  $f:\IC^{n+1}\to\IC$, $f$ 
has a 1-dimensional singular locus $\Sigma$, which is a complete 
intersection curve defined by an ideal $I$. 
Assume that $f\in I^2$. This implies that $J(f)$, the jacobian ideal 
of $f$ is in $I$ as well. Let $j(f)=\dim_{\IC}{{I}\over{J(f)}}$. This plays an important role in the theory, because it deforms flatly, if we deform $f$.

Two important examples of  non-isolated singularities are germs of type 
A$_{\infty}$ which have the normal form 
$f(z_1,\dots,z_{n+1})=\sum\limits_{i=1}^{n} {\mathop z_i^2}$ and germs of 
type
D$_{\infty}$ which have normal form 
$f(z_1,\dots,z_{n+1})=z_1z^2_2+\sum\limits_{i=3}^{n+1} {\mathop z_i^2}$. 
Note that if n=2 then D$_{\infty}$ is just a Whitney umbrella.
For A$_{\infty}$ germs $j(f)=0$ while for D$_{\infty}$ germs $j(f)=1$.

If we have a germ $f$ as in the setup for this part, then we can deform $f$ so that $F_y$ has only singularities of type A$_{\infty}$, D$_{\infty}$ and Morse singularities (${\rm A}_{1}$ points) , and the singular locus of the family of functions remains a complete intersection.  ( Cf. \cite{P1} p 83 proposition 7.18.) The existence of this deformation makes it possible to prove the following result, which unpacks the geometry of the multiplicity of the pair.

\thm 3 Suppose $f:\IC^{n+1}\to\IC$, $f$ has a 1-dimensional singular locus 
$\Sigma$, which is a complete intersection curve defined by an ideal $I$, 
$f\in I^2$ and $j(f)$ finite.
Then $$e(J(f),I)=j(f)=\#\{{\rm D}_{\infty}(F_y)\}+\#\{{\rm A}_{1}(F_y)\}.$$
where $\#\{{\rm D}_{\infty}(F_y)\}$ is the number of Whitney umbrella points that appear in a versal deformation of $f$, and $\#\{{\rm A}_{1}(F_y)\}$ is the number of Morse points.

\pf Full details are found in  \cite{G3}, but here is a sketch. Apply the multiplicity polar theorem to the deformation of Pellikaan. The induced deformation of the ideal ideal $I$ has no polar variety of the same dimension as the base of the deformation since it is a complete intersection. The polar variety of 
the deformation of $J(f)$ is empty because the deformation of $J(f)$ has $n+1$ generators. Then in \cite{G3} we show that the only points where $I_y/J(f_y)$ has support is at Whitney umbrella points, and at Morse points, and each contributes $1$ to the multiplicity  of $(I_y,J(f_y))$. This shows the equality of the ends. The second equality is due to Pellikaan, ( Cf. \cite{P1} p 87 proposition 7.20.)

Theorem 2.3 was based on a deformation of $f$, the defining equation of the hypersurface. We can also deform the parameterization of the hypersurface. This will also give us examples of surfaces whose singular locus is not a complete intersection curve, whose geometry we can describe. Now, suppose $F:\IC^2,0\to\IC^3,0$ is the germ of a finitely determined map germ. (Cf \cite {G2} for a discussion of these germs and their geometry.) They are important, because they are the type of germ which appears in generic families of map germs from $\IC^2\to\IC^3$. For this paper it suffices to know that  away from the origin, the images of these germs have only curves  of transverse double points as singular points, and that they have versal unfoldings, such that the image of the generic  member of the unfolding has only Whitney umbrella singularities, isolated traverse triple points, and curves which are transverse double points except at the Whitney umbrella points and triple points.  An unfolding of $F$ induces a deformation of $f$ the defining equation of the image of $F$. We wish to wish to apply the multiplicity polar theorem to this deformation of $f$. The singular locus of the unfolding of $F$ is defined by its conductor ideal ${\cal C}$, which is determinantal and specializes (${\cal C}(y)={\cal C}(F_y)$). For a proof of these facts see \cite{MoP}.  Since the triple points are a new type of singularity, we give them a closer look.

\lem 4 At a triple point of a multi-germ $F$ with its image defined by $f$, $J(f)={\cal C}$.

\pf We can assume $f(x,y,z)=xyz$, then $J(f)=(yz,xz,xy)={\cal C}$, since both $J(f)$ and ${\cal C}$ then define the reduced structure on the singular set of the image of $F$.

Denote a subsheaf  of  ${\cal C}$ generated by $3$ elements which define the polar variety of ${\cal C}$ of codimension 3 by  ${\cal C}_P$.

\cor 5 For a generic $y$ value in the base of the versal deformation of $F$,  at the triple points of the image of $F(y)$,  we can choose ${\cal C}_P$, so that ${\cal C}_P(y)={\cal C}(y)$.

\pf  Since the result holds for the generators of $J_z(\bar f)$ where $\bar f$ is the induced deformation of $f$, and $J_z(\bar f)$ is the ideal generated by the partial derivatives with respect to the coordinates on $\IC^3$, it holds for a good choice of ${\cal C}_P$--just consider a linear combination of three generators of ${\cal C} $ with the generators of $J_z(\bar f)$.

Applying the multiplicity polar theorem to the pair $(J_z(\bar f), {\cal C})$ viewed as ideals in $\O_{3+k}$ we obtain:

\thm 6 Suppose $f$  defines the image of $F:\IC^2,0\to\IC^3,0$, $F$ a finitely determined map germ, with a $k$-parameter versal unfolding. Then

i) $ {\rm  mult}_y \Gamma_{3}({\cal C,O}_{3+k})=\dim_{\IC}  {\cal C}(0)/{\cal C}_P(0)$

ii) $e(J(f), {\cal C}(0),0)+ {\rm  mult}_y \Gamma_{3}({\cal C})=e(J(f), {\cal C}(0),0)+\dim_{\IC}  {\cal C}(0)/{\cal C}_P(0)$

$=\dim_{\IC} {\cal C}(0)/J(f)=\#\{{\rm D}_{\infty}(f_y)\}+\#\{{\rm A}_{1}(f_y)\}.$

\pf Since $ {\cal C} $ is determinantal, it is Cohen-Macaulay; by section 4 of chapter 2 of \cite{P1}, it follows that $ {\cal C} /{\cal C}_P$ is perfect. This implies that 
$\dim_{\IC}  {\cal C}(y)/{\cal C}_P(y)$ is independent of $y$; for generic $y$ it follows that the support of $\dim_{\IC}  {\cal C}(y)/{\cal C}_P(y)$ is the polar variety $\Gamma_{3}({\cal C})$, and each point contributes $1$ to $\dim_{\IC}  {\cal C}(y)/{\cal C}_P(y)$. This is true because at points of the cosupport of ${\cal C}$, where the cosupport is smooth, in order for ${\cal C}_P(y)$ to be a reduction of ${\cal C}(y)$, they must be equal. The only points of the cosupport which are singular for generic $y$ are the triple points, and these are covered by Lemma 2.4.

Now applying the multiplicity polar theorem to the pair $(J_z(\bar f), {\cal C})$, we obtain the equality of the left and rightmost terms of ii) as the triple points do not contribute to
the cosupport of ${\cal C}(y)/J(f_y)$. The first equality follows from i) while the second follows by applying section 4 of \cite{P1} to the quotient  $({\cal C}/J_z(\bar f))$. Again this uses the fact that triple points are not in the cosupport of the quotient, and each type of point in the cosupport for generic $y$ contributes $1$ to  $\dim_{\IC}{\cal C}(y)/J(f_y)$.

The argument used to prove the second equality of ii) fills in a gap in the proof of the main theorem of \cite {Mo1}. There the reference to \cite{P1} given to establish the second equality is not really relevant to what is being proved. (However, as our argument shows, it is not difficult to make the argument needed  using other parts of \cite {P1}.)

If we look at a value of our parameter $y$ for which $F_y$ is a stable germ, then Mond has shown that the homotopy type of the image is a bouquet of spheres (\cite {Mo1} theorem 1.4 p226). In this context the image is called a disentanglement of the image of $F$. We can give a formula for the number of spheres in this bouquet using the multiplicity of the pair. Denote the number of spheres in the bouquet by $\mu(F)$. By pulling back by $F$ to $\IC^2$,  we can form the quotient  ${\cal C}(0)/J(f){\cal O}_2$.  Mond also showed that the number of Whitney umbrellas appearing in the versal deformation of  
$F$ is $\dim_{\IC}{\cal C}(0)/J(f){\cal O}_2$. Then we have:

\cor 7 Suppose $F:\IC^2,0\to\IC^3,0$ is a finitely determined germ, then 
$$\mu(F)=e({\cal C}(0), J(f),0)+\dim_{\IC}  {\cal C}(0)/{\cal C}_P(0)-\dim_{\IC}{\cal C}(0)/J(f){\cal O}_2.$$
\pf This uses Theorem 2.6 and a result of Siersma showing that the number of $\{{\rm A}_{1}(f_y)\}$ is $\mu(F)$.

In our next application, we show how the multiplicity polar theorem can be used to describe the number of singularities of a differential form on an equidimensional complex analytic germ with isolated singularities.

In a series of papers,(\cite {E-GZ1},\cite{E-GZ2},\cite{E-GZ3}, and \cite{EGZS})
 Ebeling and Gussein-Zade have been looking at the index of a differential 1-form $\omega$ defined on a singular space. In the case that
that $X$ has an isolated singularity they defined the radial index of a differential 1-form $\omega$ at a singular point $z$,
denoted ${\rm{ind}}_z(\omega)$, in \cite{E-GZ1}.
The radial index has the property that on smooth varieties it agrees with the usual notion of the index, and satisfies conservation of
number.  A differential form $\omega$ on a smooth analytic set $X$ has a singularity at a point $x$, if its restriction to the tangent space to $X$ at $x$ is zero. This means that it is a conormal to $X$, hence the matrix whose rows are $\omega$ and the rows of $DF$, $F$ the map that defines $X$ has the same rank as $DF$ at $x$. This suggests looking at the module $JM(X,\omega)$, which is the module
generated by the columns of the matrix whose last row is $\omega$, and whose first rows come from $DF$. If $X$ is an ICIS, then this module has finite colength if $\omega$ has an isolated singularity.

If $X$ is not an ICIS, then the module never has finite colength, so we look at the multiplicity of a pair.

We have that $JM(X)$ is a submodule of $F$, the free $\O_X$ module on $p_I$ generators, where   $p_I$ is the number of generators of $I$ the ideal that defines $X$. Define $H_{d-1}(X)$ as the elements of $F$ which are in the integral closure of $JM(X)$, except possibly at the origin. Since $X$ is smooth except at the origin, $JM(X)$ is free except at the origin, hence integrally closed, so $JM(X)$  and $H_{d-1}(X)$ agree except at the origin. Further, since we assume $\omega$ has an isolated singularity at $0$, we also have that  $JM(X,\omega)$ and  $H_{d-1}(X)\oplus \O_{X,0}$ agree except at the origin, so the multiplicity of the pair is well defined. By perturbing $\omega$ we can expect the singularity of $\omega$ to breakup, and we can use the multiplicity polar theorem to keep track of the contributions from the points. We say a linear form $L$ is generic for $X$ at the origin if $L$ is not a limit of tangent hyperplanes to $X$ at the origin.

 Let $\tilde\chi(M)$
 denote the
reduced Euler characteristic of $M$. Based on the results  of \cite{EGZS} we can
show: 

\thm 8 Suppose $X^d,0\subset\IC^N$ is a germ of complex analytic set with an isolated singularity at $0$. Suppose $\omega$ is a differential 1-form with
an isolated singularity at $0$. Suppose $L$ is generic linear form with respect to the singularity of $X$. Then

$${\rm{ind}}_0(\omega)=e(JM(X,\omega), H_{d-1}(X)\oplus \O_X,0)-e(JM(X,dL), H_{d-1}(X)\oplus \O_X,0)$$
$$+(-1)^{d-1}\tilde \chi(L^{-1}(t)\cap X).$$

\pf Extend the above setup to the product $X\times \IC$, replacing $\omega$ by $\omega+tdL$, $L$ generic. From conservation of number and the existence of a generic $L$ which is dealt with in \cite{EGZS}, it follows that

$${\rm{ind}}_0(\omega)=C+{\rm{ind}}_0(dL),$$
\noindent where $C$ is the number of simple singularities of $\omega+tdl$ for small $t$. Now, $H_{d-1}(X\times \IC)$, which consists of elements of the free module containing $JM(X\times\IC)$ which are integrally dependent on $JM(X\times\IC)$ except along the $t$-axis is  independent of $t$, so it has no polar curve, and the same is true for $H_{d-1}(X\times\IC)\oplus \O_{X\times\IC}$. Meanwhile, $JM(X\times\IC,\omega+tdL)$ has no polar curve either because the number of generators $JM(X\times\IC,\omega+tdL)$ is $N$ This follows because the polar curve is empty unless the dimension of the fiber of the transform of $X\times \IC$ over  the origin is at least $d+(e+1)-1=N$, since the generic rank is 
$e+1$ where $e$ is the generic rank of $JM(X)$. Since the number of generators is $N$, the fiber dimension of the transform is at most $N-1$. By the multiplicity polar theorem,
$$ e(JM(X,\omega), H_{d-1}(X)\oplus \O_X,0)=C+e(JM(X,dL), H_{d-1}(X)\oplus \O_X,0)$$
since for generic $t$ the contribution of $\omega+tdL$ to the sum at the origin is just the  value with generic $L$.

To finish the
proof, in \cite{EGZS}, in example 2.6,  it is shown that ${\rm{ind}}_0(dL)=(-1)^{d-1}\tilde \chi(L^{-1}(t)\cap X).$ 

We now give an application of the multiplicity polar theorem of the second type. Here is the set-up for this next application.

Let $(X,0)$ be a reduced equidimensional analytic germ in $(\IC^{n+k},0)$, and $(Y,0)$ a smooth sub-germ, coordinates chosen so that $Y$ is embedded in $X$ as $\IC^k\times 0$.  Let $f\colon ((\IC^{n+k},0)\to(\IC,0)$ be a map germ, $f(Y)=0$.
Assume that there is a smooth, dense, and open analytic subset $X_0$ of
$X$ such that $f|(X_0,0)$ is a submersion onto its image in $\IC$, $f$ non-constant and has
equidimensional fibers.  

 Recall from Definition~1.3.7 on p.~550 in \cite{LT} (compare with
pp.~228--29 in \cite{H-M-S}) that $(X_0,Y)$ satisfies the {\it
condition\/} {\hbox{\rm W$_f$}} {\it at\/} $0$ if there exist a
(Euclidean) neighborhood $U$ of $0$ in $X$ and a constant $C>0$ such
that, for all $y$ in $U\cap Y$ and all $x$ in $U\cap X_0$, we have
   $${\hbox{\rm dist }}\bigl( Y, T_xX(f(x))\bigr)\leq
C\,{\hbox{\rm dist }}(x,Y)\eqno(2.1.1)$$
 where  $T_xX(f(x))$ is the tangent spaces to the
indicated fibers of the restriction $f|X$.  This condition
depends only on the restriction $f|X$, and not on the embedding of $X$
into ${{\bf C}}^n$.

Say that $(X,0)$ is defined by the vanishing of $F\colon ({{\bf
C}}^{n+k},0)\to ({{\bf C}}^p,0)$.  Form the corresponding {\it augmented
Jacobian module} $JM(F;f)$: namely, first form the $p+q$ by $n$ matrix
$D(F;f)$ by augmenting the Jacobian of $F$ at the bottom with the
Jacobian of $f$; then $JM(F;f)$ is the ${{\cal O}}_X$-submodule of the
free module ${{\cal O}}_X^{p+q}$, generated by the columns of $D(F;f)$.

We have ${{\bf C}}^{n+k}={{\bf C}}^k\times \IC^n$, and form
the corresponding projections,
        $$r_k\colon{{\bf C}}^{n+k}\to \IC^k{\hbox{ and }} r_n\colon{{\bf
C}}^{n+k}\to{{\bf C}}^{n}.$$
 Form the corresponding {\it relative augmented Jacobian modules},
        $$JM(F;f)_{r_k}{\hbox{ and }} JM(F;f)_{r_n}$$
 by definition, these are the submodules of $JM(F;f)$ generated by the
partial derivatives with respect to the last $n$ variables on
${{\bf C}}^{n+k}$ and with respect to the first $k$ variables.
Finally, let ${{\bf m}_Y}$ be the ideal of $Y$ in ${{\bf C}}^{n+k}$.
 
 The relative conormal $C(X,f)$ of $f$ is just $T_{JM(F;f)}(X)$. It consists of tangent hyperplanes to the fibers of $f|(X_0)$, and their limits.
 
 Conditions like {\hbox{\rm W$_f$}}, which are defined by analytic
inequalities, often can be re-expressed algebraically in terms of
integral dependence, which we do next.

\prp 9 In the setup above, the following conditions are
equivalent:
 \smallbreak
 {\rm(i)} the pair $(X_0,Y)$ satisfies {\hbox{\rm W$_f$}} at $0;$\par
 {\rm(ii)} the module $JM(F;f)_{r_n}$ is integrally dependent on ${{\bf
m}_Y} JM(F;f)_{r_k}$;\par
 {\rm(ii$'$)} the module $JM(F;f)_{r_n}$ is integrally dependent on
${{\bf m}_Y} JM(F;f)$;\par
 {\rm(iii)} along the preimage in $C(X,f)$ of $0$, the ideal of
   $C(Y,f)\cap C(X,f)$ is integrally dependent on the ideal of the
preimage of $Y$.
\pf See the proof of Proposition 2.3 of \cite{GK2}.

Now we use the elementary version of the PSID we proved in Corollary 1.6. The invariant we will be using is $e_{\Gamma}(M,N,y)$, defined as follows:

$$e_{\Gamma}(M,N,y):=\Sum_{x\in p^{-1}(y)}e(M(y),N(y),x)+{\rm mult}_Y\Gamma(N,x)$$
where $p$ is the projection to $Y$.

It may be necessary to rechoose the representative of the polar variety at different $x$.

We will be assuming that we have a family of functions with isolated singularities, on a family of sets also with isolated singularities. With these assumptions, we have that 
the multiplicity of the pair $(JM(F_y,f_y), H_{d-1}(JM(F_y))\oplus {\cal O}_{X_y})$ is defined for $y$ sufficiently close to $0$. In applying our machinery, we need to use a perhaps smaller submodule for $N$, $H_{d-1}(JM(F)_{r_k})\oplus {\cal O}_{X}$. Generically $H_{d-1}(JM(F)_{r_k})(y)=H_{d-1}(JM(F_y))$, but for some values the module on the left may be smaller.

\thm 10  Suppose $( X^{d+k},0) \subset  (\IC^{n+k},0)$, 
$X = F^{-1}(0)$,  $F:{\IC}^{n+k}\to{\IC}^p$, $Y$ a smooth subset of $X$,
coordinates  chosen so
that ${\IC}^k \times {0} = Y$, $X$ equidimensional with equidimensional fibers, $X$ reduced, all fibers generically reduced, 
and fibers reduced over a Z-open subset of $Y$, $f:{\IC}^{n+k}\to{\IC}$, $f|X_0$ a submersion, $X_0$ a non-empty, smooth, open every where dense subset of $X$, $Z=f^{-1}(0)$.

A) Suppose $X_y$ and $Z_y$ are isolated
singularities,  suppose the singular set of $f$ is $Y$.
Suppose $e_{\Gamma}(m_YJM(F_y;f_y),H_{d-1}(JM(F)_{r_k})(y)\oplus {\cal O}_{X_y},y)$ is independent of $y$. 
 Then the union of the singular points
of $f_y$
 is $Y$, and the pair of strata
$(X-Y,Y)$ satisfies condition  {\hbox{\rm W$_f$}}.

B) Suppose $\Sigma(f)$
 is $Y$ or is empty, and  the pair $(X-Y,Y)$ satisfies {\hbox{\rm W$_f$}}.  Then   $e_{\Gamma}(m_YJM(F_y; f_y),H_{d-1}(JM(F)_{r_k})(y)\oplus {\cal O}_{X_y},y)$ is independent of $y$.
\pf The proof is similar to the proof of Theorem 5.7 of  \cite{G6}.

First we prove A). Since  {\hbox{\rm W$_f$}} holds generically, we have generically $JM(F;f)_{r_n}$ is integrally dependent on ${{\bf
m}_Y} JM(F;f)_{r_k}$; then we can apply Corollary 1.6 (the PSID), to deduce that 
${JM(F;f)_{r_n} }\subset\?{m_YJM(F;f)_{r_k}}$. This implies that $\?{JM(F;f)}\subset\?{JM(F;f)_{r_k}}$.
 Hence the union of the singular points
of $(F_y,f_y)$ which is the cosupport of $\?{JM(F;f)_{r_k}}$ is equal to the cosupport of $\?{JM(F;f)}$
which is $Y$.  Then the inclusion ${JM(F;f)_{r_n}}\subset\?{m_YJM(F;f)_{r_k}}$ implies {\hbox{\rm W$_f$}} for $(X-Y,Y)$. 

Now we prove B).  {\hbox{\rm W$_f$}} implies ${JM(F;f)_{r_n}}\subset\?{m_YJM(F;f)}$ which implies that $\?{m_YJM(F;f)}=\?{m_YJM(F;f)_{r_k}}$.
We know by \cite{H-M-S} that condition {\hbox{\rm W$_f$}} implies that the fiber dimension of the exceptional divisor of
$B_{m_Y}(C(X,f))$ over each point of $Y$ is as small as possible.
The integral closure condition $\?{m_YJM(F;f)}=\?{m_YJM(F;f)_{r_k}}$ implies that the same is true for
$B_{m_Y}T_{JM(F;f)_{r_k}}X)$. This implies that the polar of $m_YJM(F)_{r_k}$ is empty, hence by the multiplicity polar formula
the invariant $e_{\Gamma}(m_YJM(F_y; f_y),H_{d-1}(JM(F)_{r_k})(y)\oplus {\cal O}_{X_y},y)$ is independent of $y$.

At this point, one way to reduce the dependence of the invariants on the family, if $X(0)$ has a versal deformation with smooth base, is to use the versal deformation to define the correct $N$. In this case the multiplicity of the polar of $N$ becomes independent of the family, because the family embeds in a universal family. This approach is used in the Whitney case in \cite{G6}.

If we know the pair $(X-Y,Y)$ satisfies the Whitney condition, then we have another way to reduce dependence on the total space of the family. If the pair $(X-Y,Y)$ is Whitney, then it is known that for $l$ a linear form defining a hyperplane $H$ which is not a limiting tangent hyperplane to $X$ at the origin, $Y\subset H$ then $(X\cap H-Y,Y)$ also satisfies the Whitney conditions. In turn this implies that the pair $(X-Y,Y)$ satisfies the  {\hbox{\rm W$_l$}} condition.

\thm 11  Suppose $( X^{d+k},0) \subset  (\IC^{n+k},0)$, 
$X = F^{-1}(0)$,  $F:{\IC}^{n+k}\to{\IC}^p$, $Y$ a smooth subset of $X$,
coordinates  chosen so
that ${\IC}^k \times {0} = Y$, $X$ equidimensional with equidimensional fibers, $X$ reduced, all fibers generically reduced, 
and fibers reduced over a Z-open subset of $Y$, $f:{\IC}^{n+k}\to{\IC}$, $f|X_0$ a submersion, $X_0$ a non-empty, open every where dense subset of $X$,$Z=f^{-1}(0)$, and $(X-Y,Y)$ satisfies the Whitney conditions, $l$ a generic linear form in the sense of the above.

A) Suppose $X_y$ and $Z_y$ are isolated
singularities,  suppose the singular set of $f$ is $Y$.
Suppose 
$$e(m_YJM(F_y;f_y),H_{d-1}(JM(F)_{r_k})(y)\oplus {\cal O}_{X_y},y)$$
$$-e(m_YJM(F_y;l),H_{d-1}(JM(F)_{r_k})(y)\oplus {\cal O}_{X_y},y)$$
 is independent of $y$. 
 Then the union of the singular points
of $f_y$
 is $Y$, and the pair of strata
$(X-Y,Y)$ satisfies condition  {\hbox{\rm W$_f$}}.

B) Suppose $\Sigma(f)$
 is $Y$ or is empty, and  the pair $(X-Y,Y)$ satisfies {\hbox{\rm W$_f$}}.  Then   
 $$e(m_YJM(F_y;f_y),H_{d-1}(JM(F)_{r_k})(y)\oplus {\cal O}_{X_y},y)$$
 $$-e(m_YJM(F_y;l),H_{d-1}(JM(F)_{r_k})(y)\oplus {\cal O}_{X_y},y)$$
  is independent of $y$.

\pf  By the argument preceding the theorem, and B) of Theorem 2.10, we know that $e_{\Gamma}(m_YJM(F_y; l),H_{d-1}(JM(F)_{r_k})(y)\oplus {\cal O}_{X_y},y)$ is independent of $y$. Since we use the same $N$, $H_{d-1}(JM(F)_{r_k}\oplus {\cal O}_{X} $for both $f$ and $l$, the polar term coming from $N$ is the same, and cancels in the difference. Applying Theorem 2.10 implies the result.

In the event $H_{d-1}(JM(F)_{r_k})(y)=H_{d-1}(JM(F_y))$ for all $y$, we get a result independent of the total space of the family. This equality is a flatness condition on the ideal of $e\times e$ minors of $D_z(F)$, where $D_z(F)$ denotes the matrix of partials of $F$ with respect to the coordinates on $\IC^n$.

\sct References

\references

BMM
J. Brian\c con{,} P. Maisonobe and M. Merle,
 {\it Localisation de syst\`emes
diff\'erentiels, stratifications de Whitney et condition de Thom,}
 \invent 117 1994 531--50
 
 E-GZ1
 W. Ebeling, S. M. Gusein-Zade: On the index of a vector field at an 
isolated singularity. In: The Arnoldfest, edited by E. Bierstone et al., 
Fields Inst. Commun. 24, AMS, 1999, pp. 141-152

E-GZ2
 W. Ebeling, S. M. Gusein-Zade, On the index of a holomorphic 1-form on an isolated complete 
intersection singularity. Doklady Math. 64 (2001), 221-224

E-GZ3
 W. Ebeling, S. M. Gusein-Zade: Indices of 1-forms on an isolated complete 
intersection singularity. Moscow Math. J. 3, 439-455 (2003)

EGZS
 W. Ebeling, S. M. Gusein-Zade, J. Seade: Homological index for 
1-forms and a Milnor number for isolated singularities. Preprint 
math.AG/0307239

E
 D. Eisenbud.
    {\it Commutative Algebra with a View Toward Algebraic Geometry}. 
   Graduate Texts in Mathematics 150 1995, Springer-Verlag

F
 W. Fulton,
 ``Intersection Theory,''
 Ergebnisse der Mathematik und ihrer Grenzgebiete, 3. Folge
 $\cdot$ Band 2, Springer--Verlag, Berlin, 1984

 G2
   T. Gaffney,
   {\it Integral closure of modules and Whitney equisingularity,}
   \invent 107 1992 301--22

G5
 T. Gaffney,{\it The theory of integral closure of ideals and modules: Applications and new developments}, New Developments in Singularity Theory, p 379-404, ed. D. Siersma etal. Kluwer 2001.

G6 
 T. Gaffney,{\it Polar methods, invariants of pairs of modules and 
equisingularity,}  Real and Complex Singularities (Sao Carlos, 2002),
Ed. T.Gaffney and M.Ruas, Contemp. Math.,\#354, Amer. Math. Soc.,
Providence,
RI, June 2004, 113-136

G-G
 T. Gaffney and R. Gassler.
    {\it Segre numbers and hypersurface singularities}.
  Journal of Algebraic Geometry 8 1999, 695-736
  
   GK
   T. Gaffney and S. Kleiman,
   {\it Specialization of integral dependence for modules}
   \invent 137 1999 541-574

GK2
T. Gaffney and S. Kleiman, {\it $W_f$ and specialization of integral dependence for modules},  The Proceedings of the 1998 Workshop on Singularities at Sao Carlos, Real and Complex Singularities, Chapman and Hall \# 412, 1999

Gr
 Grothendieck, A.
   {\it Elements Geometrie Algebrique IV}.
 Publ. Math. IHES, 28, 1966

H-M-S
J.P.G. Henry, M. Merle, and C. Sabbah,
 {\it  Sur la condition de Thom stricte pour un morphisme analytique
complexe},
  Ann. Scient. \'Ec. Norm. Sup. {\bf17} (1984), 227--68

KT1
  S. Kleiman and A. Thorup,
  {\it A geometric theory of the Buchsbaum--Rim multiplicity,}
 \ja 167 1994 168--231

KT2
 S. Kleiman and A. Thorup,
{\it The exceptional fiber of a generalized conormal space},
``New Developments in Singularity Theory", Nato Science Series v21, 2001, Ed. Siersma, Wall, Zakalyukin, p401-404.

LT
   D. T. L\^e and B. Teissier,
 {\it  Limites de'espaces tangent en g\'eom\'etrie analytique},
 Comment. Math. Helvetici, 63 (1988), 540--78.

Mo1
D. Mond, {\it Vanishing cycles for analytic maps}, Singularity Theory and its Applications, Warwick 1989, Ed. Mond, Montaldi, SLNM \#1462 p221-234

MoP
D. Mond, R. Pellikaan, {\it Fitting ideals and multiple points of analytic mappings}, Algebraic and Analytic Geometry, Patzcuaro, 1987, Ed. E. Ramirez de Arellano, SLNM 1414, 1989

M
   D.  Mumford, {\it Algebraic geometry. I. Complex projective varieties}, reprint of the 1976 edition. Classics in Mathematics. Springer-Verlag, Berlin, 1995. x+186 pp. ISBN: 3-540-58657-1

P1
R. Pellikaan, Thesis, 1985 Utrecht

P2
R. Pellikaan, Finite determinacy of functions with non-isolated singularities,Proc. London math Soc. (3) 
57 (1988) 357-382

R
 D. Rees,
{\it Reduction of modules,}
 \mpcps 101 1987 431--49

T-1
 B. Teissier,
 Cycles \'evanescents, sections planes et conditions de Whitney,
 in {\it Singularit\'es \`a Carg\`ese,} Ast\'erisque {\bf 7--8} (1973),
 285--362.

 T-2
   B. Teissier,
   {\it Multiplicit\'es polaires, sections planes, et conditions de
 Whitney,}
   in ``Proc. La R\'abida, 1981.'' J. M. Aroca, R. Buchweitz, M. Giusti and
 M.  Merle (eds.), \splm 961 1982 314--491

\endreferences

\bye